\documentstyle[11pt]{article}
\begin{document}
\bibliographystyle{unsrt}

\def\bea*{\begin{eqnarray*}}
\def\eea*{\end{eqnarray*}}
\def\ba{\begin{array}}
\def\ea{\end{array}}
\count1=1
\def\be{\ifnum \count1=0 $$ \else \begin{equation}\fi}
\def\ee{\ifnum\count1=0 $$ \else \end{equation}\fi}
\def\ele(#1){\ifnum\count1=0 \eqno({\bf #1}) $$ \else \label{#1}\end{equation}\fi}
\def\req(#1){\ifnum\count1=0 {\bf #1}\else \ref{#1}\fi}
\def\bea(#1){\ifnum \count1=0   $$ \begin{array}{#1}
\else \begin{equation} \begin{array}{#1} \fi}
\def\eea{\ifnum \count1=0 \end{array} $$
\else  \end{array}\end{equation}\fi}
\def\elea(#1){\ifnum \count1=0 \end{array}\label{#1}\eqno({\bf #1}) $$
\else\end{array}\label{#1}\end{equation}\fi}
\def\cit(#1){
\ifnum\count1=0 {\bf #1} \cite{#1} \else 
\cite{#1}\fi}
\def\bibit(#1){\ifnum\count1=0 \bibitem{#1} [#1    ] \else \bibitem{#1}\fi}
\def\ds{\displaystyle}
\def\hb{\hfill\break}
\def\comment#1{\hb {***** {\em #1} *****}\hb }

\newcommand{\TZ}{\hbox{\bf T}}
\newcommand{\MZ}{\hbox{\bf M}}
\newcommand{\ZZ}{\hbox{\bf Z}}
\newcommand{\NZ}{\hbox{\bf N}}
\newcommand{\RZ}{\hbox{\bf R}}
\newcommand{\CZ}{\,\hbox{\bf C}}
\newcommand{\PZ}{\hbox{\bf P}}
\newcommand{\QZ}{\hbox{\bf Q}}
\newcommand{\HZ}{\hbox{\bf H}}
\newcommand{\EZ}{\hbox{\bf E}}
\newcommand{\GZ}{\,\hbox{\bf G}}
\font\germ=eufm10
\def\goth#1{\hbox{\germ #1}}
\newtheorem{pf}{Prooof}
\renewcommand{\thepf}{}
\vbox{\vspace{38mm}}
\begin{center}
{\LARGE \bf The structure of quotients of the Onsager algebra by 
closed ideals }\\[5mm]

Etsuro Date 
\footnote{Supported in part by the Grant-in-Aid for Scientific Research
(B)(2)09440014 Japanese Ministry of Education, Science, Sports and Culture.}
\\{\it Department of Mathematics \\ Graduate
School of Science \\ Osaka University \\
Osaka 560-0043, Japan \\ (e-mail:
ed@math.sci.osaka-u.ac.jp)}
\\[5mm]
Shi-shyr Roan
\footnote{Supported in part by the NSC grant of Taiwan . }
\\{\it Institute of Mathematics \\ Academia Sinica \\ 
Taipei , Taiwan \\ (e-mail: maroan@ccvax.sinica.edu.tw)} \\[5mm]
\end{center}

\begin{center}
{\it Dedicated to Professor Shunichi Tanaka on the occasion of
his sixtieth birthday}
\end{center}

\begin{abstract} 
We study the Onsager algebra from the ideal theoretic 
point of view. A complete classification of 
closed ideals and the structure of quotient
algebras are obtained. 
We also discuss the solvable algebra 
aspect of the Onsager algebra through the use of
formal Lie algebras.

\par \vspace{5mm} \noindent
 \par \noindent
1990 PACS: 02, 05 

\end{abstract}

\vfill
\eject

\section{Introduction}

By the Onsager algebra  
we mean an infinite dimensional Lie algebra
with a basis $A_m$, $G_m$~($m\in{\bf Z}$) and the
commutation relations:
\begin{eqnarray}
(a)&&\left[A_m,A_l\right]=4G_{m-l},\nonumber \\
(b)&&\left[A_m,G_l,\right]=2A_{m-l}-2A_{m+l},
\label{eq:OA}
\\
(c) &&\left[G_m,G_l\right]=0.\nonumber
\end{eqnarray}
This Lie algebra appeared in the seminal paper of
Onsager \cite{O}, in which the
free energy of the two dimensional Ising model
was computed exactly.  We shall call this
algebra the Onsager algebra following the
convention named in \cite{D90}.
In his paper Onsager exploited the
(row-to-row) transfer matrix method, by which 
the calculation 
amounts to the calculation of the largest
eigenvalue of the transfer matrix. The transfer
matrix has a form
\begin{displaymath}
(const.)e^Ae^B
\end{displaymath}
where $A$ and $B$ are matrices of degree $2^n$ ($n$ being the
number of sites on a row), which are the linear
sums of tensor products  of Pauli
matrices. By analyzing  the structure of the
algebra ( representation ) generated by $A$ and $B$
in detail, Onsager derives the algebra 
$(\ref{eq:OA})$ or its representations.
The number of sites $n$ comes in the 
representation or  
in the structure of quotient algebra as follows
\begin{displaymath}
A_{m+2n}=A_m,\quad G_{m+2n}=G_m.
\end{displaymath}
Although the structure changes slightly depending
on the parity of $n$, the resulting
representation is a direct sum of quaternions and
scalars. Utilizing this analysis Onsager computed
the largest, the second largest eigenvalues and
corresponding eigenvectors of the transfer
matrix. Afterwards, Onsager re-solves the two
dimensional Ising model by using  now famous free
fermions (Clifford algebras) with B. Kaufman five
years later \cite{K, KO}.  The method  of free fermion is a
much powerful one than that  based on the
algebra  (\ref{eq:OA}).  This might be a reason
why the Onsager algebra was forgotten for a
while. In 1980s, this algebra reappears in the
renewed interests in two-dimensional integrable 
field theories. Dolan and Grady
\cite{DG} rediscover this algebra while studying
the condition for a two-dimensional field theory
to possess a infinite number  of conservation
laws.  Subsequently this algebra again appears in
the study of integrable spin chains \cite{GR}, 
then in the superintegrable chiral Potts
model \cite{AMPTY}.
The spectrum of the superintegrable chiral Potts model is shown to obey 
certain quadratic equations by a numerical study in
\cite{AMP}. Davies \cite{D90} studied
representations of the quotient of the Onsager
algebra by an ideal generated by a linear
relation among
$A_j$s and gave an answer to this observation. 
Except for these, not so much attention was paid
on the Onsager algebra. One of the present
authors \cite{R91} found that the Onsager
algebra  has a presentation as an invariant
subalgebra of the loop algebra 
${\bf C}[t,t^{-1}]\otimes \goth{sl}_2$ by an
involution while examining the papers
\cite{D90,D91}.  A generalization of the Onsager
algebra to the case $sl_n$ was studied in
\cite{UI}, which enhances the work 
\cite{AS}. Although the paper \cite{D90} is full of inspiring content,
it contains somewhat ambiguous settings and claims from mathematical
viewpoint. This was one of the motivation for the paper \cite{R91} and
also the present one.

The present study aims to pursue the direction 
in \cite{D90,D91}, but 
stress more on the quotient algebras by
ideals rather than the representations as in
\cite{D90,R91}. In the course of our study we
focus the attention on the case when the quotient
algebras do not have the center. Such an
ideal whose quotient algebra does not
contain central elements will be called a closed
ideal in this paper.  We classify all the closed
ideals of the Onsager algebra by exploiting the
presentation of this algebra as a subalgebra of
the loop algebra mentioned above. To each closed
ideal there corresponds a ``reciprocal" polynomial
in one variable.  The structure of the
quotient algebra differs according to whether
the corresponding polynomial has
$\pm 1$ as its zeroes or not. 
If $\pm 1$ are not among zeroes of the polynomial,
the quotient is a direct  sum of the algebras 
$\left({\bf C}[u]/u^{L_j}{\bf C}[u]\right)\otimes\goth{sl}_2$ for 
positive integers $L_j$. In the case when the
polynomial is $(t \pm 1)^L$, we found by
computer experiments that the derived algebras of
the quotient algebras thus arose are identical
with a series of nilpotent Lie algebra studied
by  Santharoubane \cite{S}.  Santharoubane
obtained such a series in a project of
reducing the classification of nilpotent Lie algebras to that of 
certain ideals in the nilpotent part of Kac-Moody Lie algebras. 
The notion of roots in
nilpotent Lie algebras are used there to establish
the connection with Kac-Moody Lie algebras.

Now we outline the contents of this paper. We
start by recalling in Sect.2 some known facts
about  the Onsager 
algebra needed for our later discussion. In Sect.3,  
we establish a basic ideal-structure theorem
of the Onsager algebra, (see Theorem
\ref{th:central} in the content), which 
associates each closed ideal with a  
reciprocal polynomial (in one variable $t$). 
Using this result, the structure of quotients of
the Onsager algebra by closed ideals is reduced
to  those by ideals which are generated by 
powers of elementary reciprocal polynomials, 
$(t-a)(t-a^{-1})$ for  non-zero complex numbers
$a$. In Sect.4 and 5, we derive the structure
of those quotient algebras for $a \neq \pm1$ and 
$a = \pm 1$ respectively. Furthermore, we 
study the completion of the Onsager algebra
through these quotient algebras, which provides
a different  realization of the 
Onsager algebra in a completion of nilpotent part of $A_1^{(1)}$ in the 
principal realization.
The relation with solvable, nilpotent
algebras is also discussed in Sect.5 for the case
$a=1$. In Sect.6, we derive the
classification of finite-dimensional irreducible
representations of the Onsager algebra, a result
known in \cite{D90, R91}, from our ideal theoretic
point of view, and discuss its physical
application to the spectra of the superintegrable
chiral Potts Hamiltonian.   Finally,  we conclude
in Sect.7 with some remarks.

{\bf Notations. } To present 
our work, we prepare some notations.
In this paper, we shall use the following 
conventions:
\par \noindent
$\goth{sl}_2$ = the Lie algebra $\goth{sl}_2(\CZ)$ 
with the standard generators $e, f, h$, 
$$
[e, f] = h \ , \ \ [h, e] = 2e \ , \ \ 
[e, f] = -2f \ . 
$$  
$\theta: \goth{sl}_2 \longrightarrow \goth{sl}_2$, the
(Lie)-involution defined by $\theta(e )= f, 
\theta(f)= e, \theta(h)= -h$. \par \noindent
$L(\goth{sl}_2) = \CZ[t, t^{-1}] \otimes
\goth{sl}_2$,  the loop  algebra of $\goth{sl}_2$
with the Lie-bracket
$$
[ p (t)  x, q(t)  y] = p(t)q(t)  
[ x , y ] \ , \ \ {\rm for} \ p(t), q(t) \in
\CZ[t, t^{-1}], \ \ \ x, y \in \goth{sl}_2 \ . 
$$
$\hat{\theta}: L(\goth{sl}_2) \longrightarrow L(\goth{sl}_2)$, the
involution defined by  $\hat{\theta}(p(t) \otimes x) = 
p(t^{-1}) \otimes \theta(x)$. \par \noindent
$\goth{sl}_2 [[u]] $ = $\CZ[[u]]\otimes 
\goth{sl}_2$, the Lie algebra of formal series in
$u$
with coefficients in $\goth{sl}_2$.
\par
\noindent  For a (Lie)
ideal
$\goth{I}$ of a (non-trivial) Lie algebra 
$\goth{L}$   (over $\CZ$), we shall denote 
$$
Z(\goth{I}) := \{ x \in \goth{L} \ | \ [x, \goth{L}] \subset \goth{I}\} , 
$$
which is an  ideal of $\goth{L}$ such that 
$Z(\goth{I})/\goth{I}$ is the  center of quotient
algebra
$\goth{L}/\goth{I}$. 
A (non-trivial) ideal $\goth{I}$ is called a
closed ideal if $Z(\goth{I})= \goth{I}$,
equivalently, $\goth{L}/\goth{I}$ is a Lie
algebra with trivial center.

\section{The Onsager Algebra}
Let $A_0, A_1$ be elements of a Lie algebra and denote
$$
G_1 = \frac{1}{4} [A_1, A_0 ] \ .
$$
An infinite sequence of elements $A_m , G_m 
 ( m \in \ZZ )$ is defined by the relations
$$
A_{m-1} - A_{m+1} = \frac{1}{2} [A_m , G_1 ] \ , \ \ 
G_m = \frac{1}{4} [A_m, A_0 ]  \ .
$$
\newtheorem{theorem}{Theorem}
\begin{theorem}\label{th:Onsager}
The
following conditions are equivalent:

(I) The elements $A_0, A_1$ satisfy the
Dolan-Grady (DG) condition $\cite{DG}$:
$$
[A_1, [A_1, [A_1, A_0]]]= 16 [A_1, A_0] \ , \ \ 
[A_0, [A_0, [A_0, A_1]]]= 16 [A_0, A_1] \ .
$$

(II) The infinite sequence of elements $A_m, G_m$  $(m \in
\ZZ)$ satisfy the relation $(\ref{eq:OA})$.

\end{theorem}
{\bf Proof.} 
The proof of $(I) \Longrightarrow
(II)$ can be found in  \cite{D91, R91}, which will
be omitted here. For $(II) 
\Longrightarrow (I)$, it follows from the
relations in $(1 a)$ for
$ (m, l)  = (0,
\pm1), (1, 0)
$, and $(1 b)$ for
$(m, l)=(0, 1), (1, 1)$. 
$\Box$ \par \vspace{.2in} \noindent
By the above theorem, we introduce the following
definition of the Onsager algebra \cite{R91}.
\par \vspace{.2in} \noindent
{\bf Definition.} The Onsager algebra, denoted by
$\goth{OA}$,  is the universal Lie algebra
generated by  two elements $A_0, A_1$ with the
$DG$-condition. Equivalently,  $\goth{OA}$ is
identified 
 with the fixed Lie-subalgebra of
$L(\goth{sl}_2)$ by  the involution $\hat{\theta}$. 
The elements $A_m, G_m$ of
$\goth{OA}$ have the following 
expressions in  $L(\goth{sl}_2)$
$$
\begin{array}{lll}
A_m = 2t^m e + 2t^{-m} f \ , & 
G_m = (t^m - t^{-m}) h \ , & {\rm for} \ \ m \in \ZZ \ .
\end{array}
$$
$\Box$ \par \vspace{.2in} \noindent
We shall always make the above identification of 
$A_m, G_m$ as elements in
$L(\goth{sl}_2)$ in what follows.  For  an
element $X$ of $L(\goth{sl}_2)$, the  criterion
of $X$ in
$\goth{OA}$ is given by
$$
\begin{array}{l} 
X \in \goth{OA} \   \Longleftrightarrow  \ \ X = 
p(t)e + p(t^{-1})f +
q(t)h \ \ \ \ {\rm with} \  \ q(t)= -
q(t^{-1})  \ ,
\end{array}
$$
where $ p(t) , q(t) \in \CZ[t, t^{-1}]$. 
Note that a polynomial $q(t)$ with the above property can be written in 
the form 
$$
q(t) = q_+(t) - q_+ (t^{-1}) \ \ \ {\rm with} \ \ 
q_+(t) \in \CZ[t] \ .
$$
Then the following equalities hold:
\begin{eqnarray}
&&[ p(t)e + p(t^{-1}) f + q(t) h , \ e+f ]\nonumber\\
&&\qquad  = 
2q(t) e + 2 q(t^{-1}) f + (p(t)- p(t^{-1})) h \
, \label{eq:fm} \\ 
&&[ p(t)e + p(t^{-1}) f + q(t) h ,
\  t e+ t^{-1} f ] \nonumber\\
&& \qquad = 
 2t q(t) e  - 2t^{-1}q(t) f
+  ( t^{-1} p(t) - t p(t^{-1}) ) h \ . \nonumber 
\end{eqnarray}
It is easy to see that 
the universal enveloping algebra of $\goth{OA}$ is the 
fixed-subalgebra of the universal enveloping
algebra of $L(\goth{sl}_2)$ by
$\hat{\theta}$, hence with the inherited 
co-multiplication structure:
$$
A_m \mapsto A_m \otimes 1 + 1 \otimes A_m \ , \ \ 
G_m \mapsto G_m \otimes 1 + 1 \otimes G_m \ .
$$
In  $\goth{OA}$, there are two involutions $\iota, \sigma$ defined by 
$$
\begin{array}{lll}
\iota :  
p(t)x \mapsto p(t^{-1})x \ , & {\rm i.e.}, & 
\iota (A_m) = A_{-m} ,  \ \  \iota (G_m) =
 G_{-m} \ , \\
\sigma :  
p(t)x \mapsto p(-t)x \ , & {\rm i.e.}, & 
\sigma (A_m) =(-1)^m A_m , \ \  \sigma (G_m) =
(-1)^m G_m \ , 
\end{array}
$$
which we will use later in the paper.

\section{Ideal Structure of the  Onsager
Algebra} 
In this section, we are going to
determine the structure of closed ideals of
\goth{OA}. Let $P(t)$ be a non-trivial monic
polynomial in $\CZ[t]$. We call $P(t)$ a  
 reciprocal polynomial if $
P(t) = \pm t^d P(t^{-1})$, where $d$ is the
degree of $P(t)$. 
Then one has,  $P(t)\CZ[t, t^{-1}]=
P(t^{-1})\CZ[t, t^{-1}]$. It is easy to see that 
the zeros of $P(t)$ not equal to $\pm 1$ occur in
reciprocal pairs. In fact, $P(t)$ is a product
of the following elementary reciprocal
polynomials
$U_a(t)$ for $a \in \CZ^*$, where  $U_a(t)$ is 
defined by 
\bea(l)
U_a(t) : = \left\{ \begin{array}{ll}
t^2 -(a+a^{-1})t + 1  \  & \mbox{if } \ a^2 \neq 1 \ ,  \\
t - a   & \mbox{if } \ a = \pm 1 \ .
\end{array} \right.
\elea(Ua)
Note that $U_a(t)= U_{a^{-1}}(t)$. 
Let $P(t)$ be a reciprocal polynomial. We call an
element
$X$ of $\goth{OA}$ divisible by $P(t)$,
denoted by $P(t)| X$,  if $X =p(t)e + p(t^{-1}) f
+ q(t) h $ with $ p(t) , q(t) \in P(t)\CZ[t,
t^{-1}]$. Denote 
$$
\goth{I}_{P(t)} := \{ X \in \goth{OA} \ | \ P(t) | X \} \
.
$$
Then $\goth{I}_{P(t)}$ is an ideal in $\goth{OA}$ 
invariant under the involution $\iota$. For two  
reciprocal polynomials
$P(t)$ and $ Q(t)$, one has the relation
$$
\goth{I}_{P(t)} \bigcap \goth{I}_{Q(t)} = 
\goth{I}_{{\rm lcm}(P(t), Q(t))} \ .
$$
In particular,  $\goth{I}_{P(t)}
\subseteq \goth{I}_{Q(t)}$ if $Q(t)| P(t)$, 
hence there is the canonical  projection
$$
\goth{OA}/\goth{I}_{P(t)} \longrightarrow 
\goth{OA}/\goth{I}_{Q(t)} \ .
$$ 
The following lemma will be useful for later
purpose.
\newtheorem{lemma}{Lemma}
\begin{lemma}\label{lem:prod}
Let $P_j(t) 1 \leq j \leq J$ be pairwise
relatively prime reciprocal polynomials and  
$P(t) := \prod_{j=1}^J P_j(t)$. Then the canonical 
projections give rise to a   
Lie-isomorphism:
$$
\goth{OA}/\goth{I}_{P(t)}
\stackrel{\sim}{\longrightarrow}
\prod_{j=1}^J  \goth{OA}/\goth{I}_{P_j(t)}
$$
\end{lemma}
{\bf Proof.} The injective part is
obvious, so  we only need to show the
surjectivity of the above map. For   $X_j =
p_j(t)e+ p_j(t^{-1})f+ q_j(t)h \in \goth{OA},  1
\leq j
\leq J $, let
$N$ be a positive integer such that
$ t^N p_j(t) , t^N q_j(t)$ are all polynomials in
$t$.  By the Chinese remainder theorem, there
exist polynomials
$\tilde{p}(t), \tilde{q}(t)
\in
\CZ[t]$ such that the  following relations hold
in
$\CZ[t]$: 
$$
\tilde{p}(t) \equiv t^N p_j(t)  \ \ \ , \ \ 
\tilde{q}(t) \equiv t^N q_j(t) \ ( {\rm mod.} \
P_j(t) ) \
\ {\rm for 
\ all \ } \ j \ .
$$
Define the element $X$ of $\goth{OA}$ by 
$$
X: = p(t)e+ p(t^{-1})f +q(t)h \ , 
$$ 
where
$$
p(t) = \frac{\tilde{p}(t)}{t^N} \ , \ \ 
q(t) = \frac{1}{2}\left(\frac{\tilde{q}(t)}{t^N} -t^N
\tilde{q}(t^{-1})\right) \ \in \CZ[t, t^{-1}] \ .
$$
By $q_j(t)+ q_j(t^{-1})=0$, one can
easily see that both 
$ q(t) - q_j(t),  p(t) - p_j(t)$ are
divisible by 
$P_j(t)$ for  all $j$. Hence $X \equiv X_j$ (mod. $\goth{I}_{P_j(t)}$) 
for all $j$. 
$\Box$ \par \vspace{.2in} \noindent 
\begin{lemma}\label{lem:recip}
Let $P(t)$ be a reciprocal polynomial
and write
$$
P(t) = (t-1)^L (t+1)^K P^*(t) \ , \ \  L , 
K \geq 0 , 
\ \ P^*(\pm1) \neq 0 \ .
$$
Denote 
$\displaystyle{\widetilde{P}(t) :=
(t-1)^{2[\frac{L}{2}] } (t+1)^{2[\frac{K}{2}]
}P^*(t)}$,  here $[r]$ stands for the integral
part of a rational number 
$r$. Then 
$$
Z( \goth{I}_{P(t)} ) = \{ 
p(t)e + p(t^{-1}) f + q(t) h \in \goth{OA} \ | \ \ 
\widetilde{P}(t) | p(t) , \ \  P(t) | q(t) \ \} \ .
$$
As a consequence, 
$\goth{I}_{P(t)}$ is closed if and only if the   
zero-multiplicities of $P(t)$ at $t=
\pm1$  are even.
\end{lemma}
{\bf Proof.}  
For  
$X = p(t)e + p(t^{-1}) f + q(t) h \in \goth{OA}$,   
the criterion of $X$ in $ Z(\goth{I}_{P(t)})$
is given by the relations
$$
P(t) | \ [X, e+f], \ [X , te+t^{-1}f] \ ,
$$
which by (\ref{eq:fm}) is the same as
$$
P(t) | \ q(t) , \ (t-t^{-1})p(t), \ p(t)- p(t^{-1}) \ .
$$
Let $p(t)$ be an element of $\CZ[t, t^{-1}]$ which
satisfies 
the above  conditions. Then
$P^*(t) | p(t)$. In order to have the conclusion
of the lemma, we may assume either
$L$ or $K$ is greater than 0, say $L>0$. One
has 
\bea(l)
(t-1)^{L-1} | p(t) \ , \ \ (t-1)^L | ( p(t) - p(t^{-1})
)
\ .
\elea(t1)
Write $p(t) = (t-1)^{L-1} h(t) $. We have
$$p(t) - p(t^{-1}) = (t-1)^{L-1} 
\left(h(t) -(-t)^{-L+1} h(t^{-1})\right). $$
Then 
$(t-1) | (h(t) -(-t)^{-L+1} h(t^{-1}))$, {\it i.e.}, 
$h(1)(1 -(-1)^{-L+1} ) = 0$, which is
equivalent to the following relation
$$ 
 \left\{ \begin{array}{ll}
(t-1)^L| p(t) & \mbox{for even} \ L \ , \\
(t-1)^{L-1}| p(t) & \mbox{for odd} \ L \ .
\end{array} \right. 
$$
Therefore the relation (\req(t1)) is the same as $
(t-1)^{2[\frac{L}{2}]} | p(t)$. The same 
argument also applies  to the case $(t+1)^K$
for $K>0$. Hence we obtain the result. 
$\Box$ \par \vspace{.2in} \noindent   
\begin{lemma}\label{lem:class}
Let $I$ be an ideal in $\goth{OA}$ and  
$r(t)$ be an element of $\CZ[t, t^{-1}]$. 

(i) If $r(t)e + r(t^{-1})f $ is an element in $I$,
then  
$(p(t)-p(t^{-1}))h \in I$ for $p(t) \in 
r(t)\CZ[t, t^{-1}]$.

(ii) For a closed ideal  $I$ and an 
integer $l$ , one has
\begin{eqnarray*}
(t^jr(t)- t^{-j}r(t^{-1}))h \in I , 
( j=0, -1)&\Longrightarrow &p(t^{\pm1}) e +
p(t^{\mp1})f \in I
\\ 
 &{\rm for}& \
 p(t) \in 
r(t)\CZ[t, t^{-1}] ; 
\end{eqnarray*}
\begin{eqnarray*}
r(t)e + r(t^{-1})f \in I &\Longleftrightarrow&   r(t^{-1})e + r(t)f \in I  \\ 
&\Longleftrightarrow&    t^l r(t)e + t^{-l}
r(t^{-1})f\in I\\  &\Longleftrightarrow&  (t^j
r(t)- t^{-j} r(t^{-1}))h\in I \ , \ \ (j=0, -1) .
\end{eqnarray*}
\end{lemma}
{\bf Proof.} It is easy to see that the
equivalent relations in the second part of
$(ii)$ follows from the other ones of the lemma.
We may also assume the 
$p(t)$ in the statement with the form
$p(t)= t^m r(t)$ for $m \in \ZZ$. Write 
$
r(t) = \sum_k a_k t^k$, where $ a_k  \in \CZ , 
a_k = 0 $  for $  |k|
\gg 0$. One has 
\begin{eqnarray*}
&&2 r(t)e+ 2 r(t^{-1})f = \sum_k a_k A_k,\\
&&(r(t)- r(t^{-1}))h  = \sum_k a_k G_k , \\
&&(t^{-1}r(t)- tr(t^{-1}))h  = \sum_k a_k G_{k-1} .
\end{eqnarray*}
By the relation
$(\ref{eq:OA}) (a)
$ one obtains $(i)$. If $I$ is a closed
ideal and $\sum_m a_k G_{k+j} \in I $ for $ j=0,
-1$ by
$(\ref{eq:OA}) (a) $ we have 
$$
4\sum_m a_k G_k = \left[ \sum_m a_k A_k , A_0\right]  \ , \ \  
4\sum_m a_k G_{k-1} = \left[ \sum_m a_k A_k , A_1\right] \ , 
$$
which implies $r(t)e+ r(t^{-1})f \in I$. Using
$(i)$, one has 
$(t^{m+j}r(t)- t^{-m-j}r(t^{-1}))h \in I$ for $
j=0, \pm1 $. With the same
argument, we can also show  
$p(t^{\pm1}) e+ p(t^{\mp1})f \in I$.
$\Box$ \par \vspace{.2in} \noindent
\begin{theorem}\label{th:central}
Let $I$ be an ideal in $\goth{OA}$. 
Then $I$ is closed if  and only if $I =
\goth{I}_{P(t)}$ for a reciprocal polynomial
$P(t)$ whose zeros at
$t=\pm1$ are of even multiplicity.
\end{theorem}
{\bf Proof.} The ``if" part follows from Lemma 
\ref{lem:recip}.
Let $I$ be a closed ideal. Denote
$$
\bar{I} : = \{ r(t) \in \CZ[t, t^{-1}] \ | \ r(t)e+
r(t^{-1})f \in I \} \ .
$$ 
By Lemma \ref{lem:class} $(ii)$, $\bar{I}$ is an
ideal in
$\CZ[t, t^{-1}]$  invariant under the involution
$r(t) \mapsto r(t^{-1})$. Let $P(t)$ be the
monic polynomial which  generates the ideal
 $\bar{I}
\cap \CZ[t]$ of the polynomial ring $\CZ[t]$.
Then
$P(t)$ is a reciprocal polynomial and $
\bar{I} = P(t) \CZ[t, t^{-1}]$. By
(\ref{eq:fm}), an  element 
$q(t)h$ of $\goth{OA}$ is divisible by $P(t)$ if and
only if it  belongs  to $I$. Hence $I$
contains the ideal $\goth{I}_{P(t)}$.  We are going to
show $I=  \goth{I}_{P(t)}$, by which and Lemma
\ref{lem:recip}, the result  follows immediately.
Otherwise, there exists an element $X$ of $I
\setminus
\goth{I}_{P(t)}$ and write
$$
X = p(t)e + p(t^{-1}) e + q(t) h \ , \ q(t) =
q_+(t)-q_+(t^{-1})
\ , 
$$
where $\displaystyle{p(t) \in \CZ[t, t^{-1}] \ , \ q_+(t)\in
\CZ[t]}$.
Note that $q(t)$ is not divisible by $P(t)$,
neither for $p(t)$ by the first equality in
(\ref{eq:fm}). We may assume 
the polynomial
$q_+(t)$ with the degree less than that of
$P(t)$ and   
$q_+(0)=0$. Let
$\tilde{X}$ be a such element $X$ with the 
degree of
$q_+(t)$ being the maximal one. By
(\ref{eq:fm}), we have the following expressions
of elements in $I$
$$
\begin{array}{ll}
[\tilde{X}, e+f] &= 2q(t)e+ 2q(t^{-1})f+ k(t)h \ , \\ 

[[\tilde{X}, e+f], te+t^{-1}f]&= 2tk(t)e- 2t^{-1}k(t)f +
2(t^{-1} q(t) -tq(t^{-1}) h .
\end{array}
$$
Note that both $q(t)$ and $ k(t) $ are not
divisible by
$P(t)$. One can write $
t^{-1}q(t) -tq(t^{-1})= \widetilde{q}(t) -
\widetilde{q}(t^{-1}) $, where  
$\widetilde{q}(t) := tq_+(t) + t^{-1}q_+(t) \in
\CZ[t]$.  The degree of the polynomial 
$\widetilde{q}(t)$ is greater than that of
$q_+(t)$, so is $\widetilde{q}(t)- \widetilde{q}(0)$.
By the definition of
$\tilde{X}$, $ t^{-1}q(t) -tq(t^{-1})$ is divisible by
$P(t)$,  which implies $P(t) | tk(t)$. By $P(0) \not=
0$,
$k(t)$ is divisible by $P(t)$, which contradicts 
our previous statement for $k(t)$. 
$\Box$ \par \vspace{.2in} \noindent
{\bf Remark.} The ideal $\goth{I}_{P(t)}$ in the
above theorem is characterized as the minimal
closed ideal of $\goth{OA}$ containing $P(t)e+
P(t^{-1})f$. We shall call $P(t)$ the generating
polynomial of the  closed ideal
$\goth{I}_{P(t)}$.  
$\Box$ \par \vspace{.2in} \noindent  

\noindent{\bf Remark.} For affine Lie algebras classification of ideals
is known (cf. Theorem 4 of \cite{Mo} and Lemma 8.6 of \cite{Ka}).
$\Box$ \par \vspace{.2in} \noindent  

An important class of closed ideals arises from 
representations of $\goth{OA}$, by which 
we shall always mean  finite dimensional
Lie-algebra representations in this paper.
The 
kernel of a representation $\rho: \goth{OA}
\longrightarrow \goth{gl}(V)$ is always an ideal
of
$\goth{OA}$, which will be denoted by ${\rm
Ker}(\rho)$. As the
$DG$-condition is  unchanged by adding constants
on
$
\rho(A_0) $ and $\rho(A_1)$,  one may assume the
representation takes  values in $\goth{sl}(V)$. For an
irreducible representation
$\rho$ of 
$\goth{OA}$ in $\goth{sl}(V)$,  by Schur's Lemma
${\rm Ker}(\rho)$ is a closed ideal, hence
invariant under the involution $\iota$ by Theorem 2. 
A representation $\rho$ of $\goth{OA}$ on 
a vector 
space $V$ is called  Hermitian if both
$\rho(A_0)$ and $\rho(A_1)$ are Hermitian operators on $V$, which
is equivalent to  the Hermitian
property of $\rho(A_j)$ for all $j$. Note that in
this  situation, the  operators $\sqrt{-1}\rho(A_0)$ and 
$\sqrt{-1}\rho(A_0)$ give rise to a
representation of $\goth{OA}$ into
$\goth{su}(V)$. Hence every Hermitian
representation of
$\goth{OA}$  is completely reducible.  For 
an irreducible Hermitian representation $\rho$ of 
$\goth{OA}$ in $\goth{sl}(V)$, ${\rm Ker}(\rho)$ is invariant 
under complex conjugation,  equivalently the 
generating polynomial $P(t)$ of ${\rm Ker}(\rho)$ 
has real coefficients.     

By Theorem \ref{th:central} and Lemma
\ref{lem:prod}, the understanding of the complete
structure of quotient algebras of
$\goth{OA}$ by closed ideals $I$ is reduced to the
case
$I=  \goth{I}_{U^L_a(t)}$ for $a \in \CZ^*, L \in
\ZZ_{\geq 0}$, where $U_a(t)$ is the elementary 
reciprocal polynomial (\req(Ua)). For $a \in
\CZ^*$, we are going to define a completion $\widehat{\goth{OA}_a} $ of
$\goth{OA}$ as follows. Set
$
\goth{OA}_{a,L}:=\goth{OA}/\goth{I}_{U^L_a(t)}$ 
and we have the canonical projections
$$
\begin{array}{ll}
\pi_{a, L} &: \goth{OA} \longrightarrow
\goth{OA}_{a,L} \\
\pi_{a, KL} &: \goth{OA}_{a,L}
\longrightarrow \goth{OA}_{a,K} \ , \ L \geq K
\geq 0 \ .
\end{array}
$$
The projective system of Lie
algebras $(\goth{OA}_{a,L}, 
\pi_{a, KL})_{L, K \in \ZZ_{\geq 0}}$ gives rise to
a limit, which is a Lie algebra 
denoted by 
$$
\widehat{\goth{OA}}_a : = \lim_{\leftarrow}
\goth{OA}_{a, L} \ . 
$$
For $L \geq 0$, there is a canonical
morphism $
\psi_{a,L}: \widehat{\goth{OA}_a} \longrightarrow 
\goth{OA}_{a,L}$. We denote its kernel 
by 
$\widehat{\goth{OA}}_a^L := {\rm
Ker}(\psi_{a,L})$. The ideals
$\widehat{\goth{OA}}_a^L$ form a filtration of
in $\widehat{\goth{OA}}_a$
$$
\widehat{\goth{OA}}_a = \widehat{\goth{OA}}_a^0
\supset 
\widehat{\goth{OA}}_a^1 \supset \cdots \supset
\widehat{\goth{OA}}_a^L 
\supset \cdots \ \ \ \ ; \ \ \ \  \ \
\widehat{\goth{OA}}_a/\widehat{\goth{OA}}_a^L \
\simeq \ \goth{OA}_{a, L}.
$$
The family of morphisms $\{\psi_{a,L}\}$ gives
rise to a morphism from $\goth{OA}$ into $
\widehat{\goth{OA}}_a$
$$
\pi_a: \goth{OA} \longrightarrow
\widehat{\goth{OA}}_a \ ,
$$
with $\psi_{a,L} \pi_a = \pi_{a,L}$. In the next
two sections, we are going to determine the
structure of $\goth{OA}_{a,
L}$ and $\widehat{\goth{OA}}_a$. For $L=1$, 
$\goth{OA}_{a,1}$
can be realized in $\goth{sl}_2$ through the
evaluation morphism of $\goth{OA}$
$$
ev_a : \goth{OA} \longrightarrow
\goth{sl}_2 \ , \ \ \ \ p(t)  x \mapsto p(a)x \ .
$$
In fact, we have the following result.
\begin{lemma}\label{lem:ker} 
${\rm Ker}(ev_a) =
\goth{I}_{U_a(t)}$ and the map $ev_a$ induces the
isomorphism: 
$$
\goth{OA}_{a, 1} \ \simeq \ 
\left\{ \begin{array}{ll} \CZ (e+f) & \mbox{if} \ a =
\pm 1 \\
\goth{sl}_2 & \mbox{otherwise} \ .
\end{array} \right.
$$
\end{lemma}
{\bf Proof.} By $U_a(a^{\pm 1})=0$, ${\rm
Ker}(ev_a)
\supseteq \goth{I}_{U_a(t)}$. Let $X$ be an
element of ${\rm Ker}(ev_a)$ with the expression 
$X= p(t)e
+ p(t^{-1}) f + q(t) h $. Then 
$p(a^{\pm 1})  = q(a)=0$. By 
$q(t^{-1})= -q(t)$,  $p(t)$ and $
q(t)$ are divisible by $ U_a(t) $, therefore $X
\in \goth{I}_{U_a(t)}$. When 
$a = \pm 1$, it is easy to see that  
${\rm Im} (ev_a)$ is the one-dimensional space
generated by
$e+f$. For $a \neq \pm 1$, ${\rm Im} (ev_a)$ is a Lie
subalgebra of $\goth{sl}_2$ containing 
$e+f, ae+a^{-1}f$.  
This implies $e, f \in  {\rm
Im}(ev_a)$, hence
${\rm Im}(ev_a)=
\goth{sl}_2 $. In fact,  the basis of
$\goth{sl}_2$ has the following expression in
term of elements of  ${\rm Im}(ev_a)$,
$\overline{A}_m:= ev_a(A_m)$:
$$
e= \frac{\overline{A}_1 -a^{-1}\overline{A}_0}{2
(a - a^{-1})}
\ ,
\
\  f = \frac{\overline{A}_1 -a\overline{A}_0}{-
2 (a-a^{-1})} \ , \ \ h = \frac{\overline{G}_1
}{a-a^{-1}} \ ,
$$
or  
$$
e= \frac{a\overline{A}_m -\overline{A}_{m-1}}{2
a^m(a-a^{-1})}
\ , \ \  f = \frac{a^{-1}\overline{A}_m
-\overline{A}_{m-1}}{- 2 a^{-m}(a-a^{-1})} \ ,
\ \  h = \frac{\overline{G}_m}{a^m - a^{-m}} \ .
$$
$\Box$ \par \vspace{.2in} \noindent

\section{Structure of the quotient by ideal
generated by $U_a(t)^L, a \neq \pm 1$ }
In the following discussion, we use the following notation
for the shifted factorial:
\begin{displaymath}
x^{(0)}=1, \quad x^{(n)}=x(x-1)\cdots (x-n+1),
\ n\in \ZZ_{>0} .
\end{displaymath}
Then 
\begin{displaymath}
{x\choose n}={x^{(n)}\over k!}.
\end{displaymath}
For $a \in \CZ^*$, we define an injective (Lie)
morphism from $\goth{OA}$ into $\goth{sl}_2[[ u]]$
by the Taylor series expansion of a function
around $t=a$ with the variable 
$u=t-a$,
\begin{equation}
s_a : \goth{OA} \longrightarrow \goth{sl}_2[[
u]] \ , \ \ \ p(t)x \mapsto \sum_{j \geq 0}
 \frac{d^j p}{dt^j}(a)
\frac{u^j}{j!} x \ , \label{eq:ta}
\end{equation}
where $p(x) \in \CZ[t,t^{-1}], x \in \goth{sl}_2$.
For a positive integer $L$, $\goth{sl}_2[[ u
]]/u^L \goth{sl}_2[[u ]] $ is isomorphic to
the Lie algebra $(\CZ[u]/u^L\CZ[u])
\otimes \goth{sl}_2$. By  composing 
with the natural projection of 
$\goth{sl}_2[[ u ]]$ to $\goth{sl}_2[[ u ]]/
u^L \goth{sl}_2[[u ]] $, $s_a$ gives rise to the
morphism
$$
s_{a,L} : \goth{OA} \longrightarrow \goth{sl}_2[[
u]]/
u^L \goth{sl}_2[[u ]] \ .
$$
Note that $s_{a,1}$ is equivalent to the
evaluation morphism $ev_a$. We are going to
extend the result of Lemma \ref{lem:ker} to an
arbitrary positive integer $L$. First we need the
following lemma.
\begin{lemma}\label{lem:invf}
 Let $f(t)$ be an
analytic function on 
$\CZ^*$. Denote $ \breve{f}(t)= f(t^{-1})$.
For 
$a \in \CZ^*$, the zero multiplicities of $f(t)$
at
$a^{-1}$ and  $\breve{f}(t)$ at $a$ are the same.
\end{lemma}

\newtheorem{proposition}{Proposition}
\begin{proposition}\label{prop:saL}
 For $L \in
\ZZ_{>0}$, Ker$(s_{a,L}) = \goth{ I}_{U_a^L}
$. The morphism $s_{a,L}$ is surjective if and
only if $a
\neq
\pm 1$.
\end{proposition}
{\bf Proof.} As before we denote $\breve{f}(t)
=f(t^{-1})$ for $f(t) \in \CZ[t,t^{-1}]$. An element 
of $\goth{OA}$ is expressed by
$
p(t) e + \breve{p} (t)f + q(t) h $ with $ 
\ p(t), q(t) \in \CZ[t,t^{-1}] $ and $  q(t)=
-\breve{q}(t) $. 
By Lemma \ref{lem:invf}, we have 
\begin{eqnarray*}
p(t)e + \breve{p}(t)f + q(t) h \in \goth{I}_{U_a^L(t)} 
&\Longleftrightarrow& p(t), q(t) \in (t-a^{\pm1})^L\CZ[[t-a^{\pm1}]] \\
&\Longleftrightarrow& p(t), \breve{p}(t), q(t) \in (t-a)^L\CZ[[t-a]] \\
&\Longleftrightarrow& p(t)e + \breve{p}(t)f + q(t) h\in  {\rm Ker}(s_{a,L}) .
\end{eqnarray*}
Therefore $\goth{ I}_{U_a^L} = {\rm Ker}(s_{a,L})$. 
For $a \neq \pm 1$,  $U_a^L(t)$
is a polynomial  of degree $2L$. One has the 
relation
\begin{displaymath}
\pmatrix{s_{a,L}(A_0) \cr s_{a,L}(A_1) \cr \vdots 
\cr \vdots \cr s_{a,L}(A_{2L-2}) \cr
s_{a,L}(A_{2L-1})}=2(C_+,C_-)
\pmatrix{e\cr\vdots\cr \frac{u^{L-1}}{(L-1)!}e
\cr f\cr\vdots\cr
\frac{u^{L-1}}{(L-1)!}f}.
\end{displaymath}
where $C_{\pm}$ are the $2L \times L$-matrices 
defined by
$$
\begin{array}{l}
 C_{\pm}  =  \left(
\begin{array}{ccc}
  a_{0,0}&  \cdots& a_{0,L-1}  
\\
a_{\pm1,0}&  \cdots& a_{\pm 1,L-1} 
\\
\vdots &  \cdots&  \vdots
\\
\vdots &  \cdots&   \vdots
\\
a_{\pm(2L-2),0}&  \cdots&a_{\pm(2L-2),L-1}  \\ 
a_{\pm(2L-1),0}&  \cdots&a_{\pm(2L-1),L-1}      
\end{array} \right)  \ , \ \ \ 
a_{ m, j}:=m^{(j)}a^{m-j} \ .
\end{array}
$$
Claim: The 
square  matrix $(C_+, C_-)$ is invertible. 
Otherwise there exists a non-trivial linear
relation  
$\displaystyle{\sum_{m=0}^{2L-1}\alpha_m a_{\pm m, j}=0,
( j=0, \ldots,
L-1)}$ for $\alpha_m \in \CZ$ not all zeros.
Then $\displaystyle{\sum_{m=0}^{2L-1}\alpha_m t^m}$ is a 
non-trivial polynomial, denoted
by
$Q(t)$. By Lemma  \ref{lem:invf}, the
zero-multiplicties of the polynomial $Q(t)$ at
$t=a^{\pm 1}$ are both at least
$L$, which contradicts the fact that the degree of $Q(t) < 2L$.
By inverting the square matrix  $(C_+, C_-)$, the
elements  
$\frac{u^j}{j!} e$ and $\frac{u^j}{j!} f$ 
for $ j=0, \ldots, L-1$, are in the Lie
subalgebra
${\rm Im }(s_{a,L})$ of $\displaystyle{\goth{sl}_2[[u]]/
u^L \goth{sl}_2[[u]]}$. Hence it follows
immediately that  the surjectivity of $s_{a,L}$
for $a \neq \pm1$. For $a = \pm1$,  we
consider the commutative  diagram
$$
\begin{array}{ccc}
\goth{OA} 
 &\stackrel{s_{a,L}}{\longrightarrow} &
\goth{sl}_2 [[u]]/u^L \goth{sl}_2 [[u]] \\
 \| & & \downarrow \\
\goth{OA} 
 &\stackrel{s_{a,1}}{\longrightarrow} &
\goth{sl}_2 \ ,
\end{array}
$$
where the right vertical morphism is the canonical
projection. By Lemma \ref{lem:ker},  the
image of
$s_{\pm1 ,1}$ is an  one-dimensional subspace of 
$\goth{sl}_2$ and it
implies the non-surjectivity of
$s_{\pm1 ,L}$.
$\Box$ \par \vspace{.2in} \noindent
As a corollary of Proposition \ref{prop:saL}, for 
$a \neq \pm 1$ 
 the morphism 
$s_{a,L}$ induces an
isomorphism,
$$
\goth{ OA}_{a, L}   \simeq \goth{sl}_2
[[u]]/u^L \goth{sl}_2 [[u]] \ , \ \ L \ \in
\ZZ_{>0} \ .
$$
By the definition of the formal algebra 
$\widehat{\goth{ OA}}_a $, one obtains the
following result.
\begin{theorem}\label{th:OAaL} For $a \in \CZ^{*}$
and
$ a
\neq 
\pm 1
$,  the morphism $s_a$ of (\ref{eq:ta}) gives rise to
the following isomorphisms,
$$
\begin{array}{ll}
\widehat{\goth{ OA}}_a \ \ \simeq  \ \ 
\goth{sl}_2 [[ u ]]
\ , 
\ \ 
& 
\widehat{\goth{ OA}}^L_a \ \ \simeq \ \  u^L
\goth{sl}_2 [[ u ]] \ ,
\end{array}
$$
and
$$
\goth{ OA}_{a,L}   \ 
 \simeq \ (\CZ[u]/u^L\CZ[u]) \otimes \goth{sl}_2
\ .
$$ 
\end{theorem}
{\bf Remark.} For $a = \sqrt{-1}$, one has 
$U_a(t)= t^2+1$.  The above structure of $\goth{
OA}_{\sqrt{-1},L} $ for
$L=2$ appeared in \cite{R91}.
$\Box$ \par \vspace{.2in} \noindent
Through the morphism $\pi_a$, one can regard
$\goth{OA}$ as a subalgebra of $\goth{sl}_2[[u]]$,
which is identified with $\widehat{\goth{OA}}_a$
by Theorem \ref{th:OAaL}.  Then the expressions of
$A_m , G_m$ in
$\goth{sl}_2[[u]]$ are given by  
\begin{eqnarray*}
&&A_m =2\sum_{j \geq 0} \left(m^{(j)}a^{m-j} 
\frac{u^j e}{j!} 
+ (-m)^{(j)}a^{-m-j}\frac{u^j f }{j!} \right) 
\ ,  \\
&&G_m =\sum_{j \geq 0} \left(m^{(j)}a^{m-j}
- (-m)^{(j)}a^{-m-j} \right) \frac{u^j h }{j!} \ .
\end{eqnarray*}
With the identification $\goth{OA}_{a,L} $  with
$\goth{sl}_2 [[u]]/u^L \goth{sl}_2 [[u]] $, one
has the expression of elements of the Onsager
algebra
\begin{eqnarray*}
&&A_m  = 2 \sum_{j=0}^{L-1} (m^{(j)}a^{m-j} 
e_j
+ (-m)^{(j)}a^{-m-j}f_j ) 
\ , \\
&&G_m  = \sum_{j=0}^{L-1} (m^{(j)}a^{m-j}
- (-m)^{(j)}a^{-m-j} ) h_j  \ ,
\end{eqnarray*}
where $e_j, f_j, h_j$ are the elements in 
$\goth{sl}_2 [[u]]/u^L\goth{sl}_2 [[u]]$
represented by 
$\frac{u^j}{j!}e, \frac{u^j}{j!}f,
\frac{u^j}{j!}h$ respectively. In fact, one
can start with  the above relations. Using the
Onsager's relation (\ref{eq:OA}) of
$A_m, G_m$ , one obtains the
relations of $e_j, f_j, h_j$ by the
following formula of shift factorials 
$$
(x+y)^{(m)} = \sum_{k = 0}^m 
{m\choose k}  x^{(k)}y^{(m-k)} \ .
$$

\section{Structure of the quotient by ideal
generated by $(t \pm 1)^L$ }
In this section, we shall discuss the structure
of 
$\goth{OA}_{a, L}$ for $a = \pm 1$. 
As the
involution $\sigma$ of $\goth{OA}$, which sends
$t$ to $-t$, induces an isomorphism
$$
\goth{OA}_{1,L} \ \ \simeq \ \ 
\goth{OA}_{-1, L} \ ,
$$
we only need to consider the case $a = 1$. For
the 
simplicity of notations, we shall omit the index 
$a=1$ in this section and denote the algebras as
$$
\goth{OA}_L = \goth{OA}_{1,L} \ , \ \ 
\widehat{\goth{OA}} = \widehat{\goth{OA}}_1 \ , \
\ 
\widehat{\goth{OA}}^L = \widehat{\goth{OA}}^L_1 \
\ ,
$$
and the morphisms as
$$
 \pi_L : \goth{OA}
\longrightarrow
\goth{OA}_L  
\ , \ \ \ 
\pi : \goth{OA}
\longrightarrow
\widehat{\goth{OA}} \ , \ \ \ 
\psi_L : \widehat{\goth{OA}} \longrightarrow
\goth{OA}_L 
\ .
$$
For  $X \in \goth{OA}$, 
the element $\pi (X)$ in 
$\widehat{\goth{OA}}$ will be  
denoted by  $X$ again later in the discussion. 
\begin{lemma}\label{lem:basis}
There are unique elements
$\underline{X}_k, \underline{Y}_k \  (0
\leq k < L)$ in $\goth{OA}_L$ such that 
\begin{displaymath}
\pi_L(A_m) = \sum_{k=0}^{L-1} {m\choose k}\underline{X}_k \
, \ \ \ 
\pi_L(G_m ) = \sum_{k=0}^{L-1}
(-1)^k{m\choose k}\underline{Y}_k \ \ for \ m \in
\ZZ \ . 
\end{displaymath}
Furthermore we have  $\displaystyle{\underline{Y}_k :={(-1)^k\over 4}
[\underline{X}_k, \underline{X}_0  ]}$.
\end{lemma}
{\bf Proof.} The above relations for 
$0 \leq m <L$ imply the uniqueness of 
$\underline{X}_k, 
\underline{Y}_k$. In fact, one can 
define $\underline{X}_k$ in terms of 
$\pi_L(A_m), (
0 \leq m <L)$ through these relations and 
set 
$\displaystyle{\underline{Y}_k :=  {1\over 4}
[\underline{X}_k, \underline{X}_0  ]}$.  By the 
definition 
of the ideal $\goth{I}_{(t-1)^L}$, we have
$$
\sum_{k=0}^{L-1} (-1)^k { L \choose k } \pi_L(A_{k+m}) = 0 \ ,
\ \ m \in \ZZ \ . 
$$
By using $(\ref{eq:OA})(a)$, it follows the 
results.
$\Box$ \par \vspace{.2in} \noindent
The above elements $\underline{X}_k,
\underline{Y}_k $ of $\goth{OA}_L$ 
for  $L > 0$ give rise to an infinite sequence of 
elements in $\widehat{\goth{OA}}$, 
$X_k, Y_k ,  ( k
\in
\ZZ_{\geq 0})$ such that the
following relations hold in $\widehat{\goth{OA}}$
\begin{equation}
A_m (= \pi (A_m)) =  \sum_{k\geq 0} {m\choose k}X_k \ , \ \ 
\ \ G_m (= \pi (G_m)) = 
\sum_{k\geq 0} (-1)^k {m\choose k}Y_k \
.\label{eq:nil}
\end{equation}
In $\goth{OA}_L$, we have
$$
\psi_L( X_k) = \left\{ \begin{array}{ll}
\underline{X}_k &
{\rm for} \  0 \leq k <L  \\
 0 & {\rm otherwise} \ . 
\end{array}
\right. 
$$
For the
simplicity of notations, the element
$\psi_L(X)$ of $\goth{OA}_L$ for $X \in
\widehat{\goth{OA}}$ will sometimes again be
denoted by
$X$ in the later discussions if no
confusion could arise. 

Recall that
the Stirling numbers of the first kind 
$$s^n_k \ \ (n, k \in\ZZ_{\geq 0})$$
are integers such that the following
relations hold for the shifted factorial
$x^{(n)}$ and the $k$-th power $x^k$
$$
x^{(n)}= \sum_{k \geq 0} x^k s^n_k.
$$
Then one has
\begin{displaymath}
{x\choose n}=\sum_{k\geq 0}{1\over k!}s_k^n x^k.
\end{displaymath}
Substituting the relations (\ref{eq:nil}) into the defining relation
of the Onsager algebra (\ref{eq:OA}) and
comparing the coefficients of $m^al^b$, we have
\begin{eqnarray*}
&&\sum_{n,k\geq 0}{1\over n!k!}s_a^ns_b^k[X_n,X_k]
=4(-1)^b{a+b\choose b}
\sum_{k\geq 0}{(-1)^k\over k!}s_{a+b}^kY_k,\cr
&&\sum_{n,k\geq 0}{(-1)^n\over
n!k!}s_a^ns_b^k[Y_n,X_k]= 2(1-(-1)^a){a+b\choose
b}
\sum_{k\geq 0}{1\over k!}s_{a+b}^kX_k,\cr
&&[Y_n,Y_k]=0.
\end{eqnarray*}

The Stirling numbers of the second kind $S^n_k (n, k \in
\ZZ_{\geq 0})$ satisfy the relation 
$$
x^n=\sum_{k \geq 0} x^{(k)}S^n_k .
$$ 
(for the general discussion of Stirling's numbers, see, e.g., 
references \cite{C, Ri} ). 
Note that $s^n_n = S^n_n = 1$, and $s^n_k =S^n_k = 0$ for $ k>n $. 
As a direct consequence of the definition
the matrices of infinite size 
$\displaystyle{\left(s_k^n\right)_{n,k\geq 0}}$
and 
$\displaystyle{\left(S_k^n\right)_{n,k\geq 0}}$
are inverse to each other
\begin{equation}
\sum_{k}  s^a_kS^k_b = \delta^a_b.\label{eq:sSinv}
\end{equation}
For Stirling numbers we have the following
identities:
\begin{lemma}\label{lem:lah}
\begin{eqnarray*}
&&\sum_{k} (-1)^k s^a_kS^k_b = (-1)^a {a!\over b!}{ a-1 \choose b-1 },
\quad {\rm \cite{Ri} \ (p.44) } \cr
&&{ j+k \choose k } S^j_a=\sum_{l} { a+l \choose l }s^l_kS^{j+k}_{a+l}.
\quad {\rm \cite{Ri2} \ (p.204) }
\end{eqnarray*}
\end{lemma}
Note that the numbers defined by the first relation are known as the Lah 
numbers.
Using (\ref{eq:sSinv}) and the identities in 
Lemma
\ref{lem:lah}, we get
\begin{proposition}
\begin{eqnarray}
&[X_n, X_k ]&= 4 (-1)^n \sum_{a \geq 0} 
{ a+k-1 \choose a } Y_{a+n+k}  \nonumber \\
&&= 4 (-1)^n \left(Y_{n+k} + \sum_{a > 0} 
{ a+k-1 \choose a } Y_{a+n+k} \right) \ ,\nonumber  \\
&[Y_n, X_k ]&= 2 \left((-1)^n X_{n+k} - \sum_{a \geq 0}
(-1)^a { a+n-1 \choose a }  X_{a+n+k} \right) \label{eq:XY} \\
&&= 2 \left(\left((-1)^n-1\right) X_{n+k} - \sum_{a >0}
(-1)^a { a+n-1 \choose a }  X_{a+n+k} \right) \ , \nonumber \\
&[Y_n, Y_k ] &= 0 \ .\nonumber 
\end{eqnarray}
\end{proposition}
With the infinite formal sum,
$\widehat{\goth{OA}}$ is a  formal Lie algebra
with generarors $X_k, Y_k$; while 
$\widehat{\goth{OA}}^L$ is  the ideal generated by
$X_k, Y_k$ with $k \geq L$. 
Hence for $L\geq
0$, 
$\widehat{\goth{OA}}^L/\widehat{\goth{OA}}^{L+1}$ 
is abelian and
$\widehat{\goth{OA}}/\widehat{\goth{OA}}^L$ is a
finite-dimensional solvable Lie algebra.   The
elements
$Y_k$ are not independent in
$\widehat{\goth{OA}}$.  We are now going to
describe the relations among $Y_k$s. Since
the commutator is skew symmetric $[X_n, X_k ]=
-[X_k, X_n ]$, from  the first relation
of  (\ref{eq:XY}) we have 
\begin{equation}
\sum_{a \geq 0} \left({ a+n-k-1 \choose a }   + (-1)^n 
{ a+k-1 \choose a } \right) Y_{a+n} = 0 \ , \ \ \ 0 \leq k \leq n \ .
\label{eq:Ynk}
\end{equation}
Since one easily sees that
 $(\ref{eq:Ynk})_{n+1,k+1}= (\ref{eq:Ynk})_{n,k+1}-
(\ref{eq:Ynk})_{n,k}$, and  $(\ref{eq:Ynk})_{n, n}=
 (\ref{eq:Ynk})_{n, 0} $, 
the relation (\ref{eq:Ynk})
is  reduced to the following one
\begin{equation}
\left(1 + 
(-1)^n \right) Y_n  + 
\sum_{a \geq 1}
{ a+n-1 \choose a }   Y_{a+n} = 0 \ , \ {\rm for } \ n \geq 0\label{eq:constraint}
\ .
\end{equation}
In particular, we have 
\begin{equation}
Y_0 = 0 \ , \ \ \
2 Y_{2n}  + \sum_{a \geq 1} { a+2n-1 \choose 2n-1 }   Y_{a+2n} = 0 \ , \ 
{\rm for } \ n \geq 1\ .
\end{equation}
In order to solve this constraint, we first
consider  their relations in 
$\goth{OA}_L$.
In this situation, we easily see that $Y_L=0$ if 
$L$ is odd.
For simplicity, we assume at this moment that $L$
is even. From (\ref{eq:constraint}), we see that
$Y_{2j}$ is a linear combination of $Y_{2k+1}$,
$j\leq k\leq [{L-1\over 2}]$ and write
\begin{equation}
Y_{2j}=\sum_{k=j}^{[{L-1\over 2}]}
\gamma_{jk}Y_{2k+1}.\label{eq:relY}
\end{equation}
After some calculation, we find that $\gamma_{jk}$ has the
following form:
\begin{equation}
\gamma_{jk}={(-1)^{k-j-1}\over 2(k-j+1)}{2k\choose 2j-1}\alpha_{k-j}.
\label{eq:form}
\end{equation}
Substituting (\ref{eq:form}) into (\ref{eq:relY}), we have
\begin{equation}
\sum_{j=0}^k{(-1)^j\over 2j+2}{2k+2\choose 2j+1}\alpha_j=1,\quad 0\leq j\leq 
[{L-1\over 2}].\label{eq:alpha}
\end{equation}
Note that this relation does not contain $L$.
In order to solve (\ref{eq:alpha}) we employ the techniques from 
inversion relations (see
\cite{Ri2}, p.109).
\begin{lemma}\label{lem:inversion}
If a sequence $\{a_k\}$ is expressed in terms of 
another sequence
$\{b_k\}$ as 
\begin{displaymath}
a_{2n+1}=\sum_{k=0}^n{2n+2\choose 2k+1}b_{2n+1-2k},
\end{displaymath}
then we have the inversion formula
\begin{displaymath}
(2n+2)b_{2n+1}=\sum_{k=0}^n{2n+2\choose 2k}d_{2k}a_{2n+1-2k},
\end{displaymath}
where $d_j$ are defined by the following expansion:
\begin{displaymath}
{2x\over e^x-e^{-x}}=\sum_{j=0}^\infty {d_{2j}x^{2j}\over (2n)!}.
\end{displaymath}
\end{lemma}
The number $d_{2j}$ is related to the 
Bernoulli number $B_j$ by the
relation 
$$d_{2j}=(-1)^j(2^{2j}-2)B_{2j},\quad j\geq 1.$$
Recall that the Bernoulli numbers $B_j$ are 
defined by the expansion
$${x\over e^x-1}= \sum_{j=0}^{\infty}{b_j\over j!}x^j~\hbox{~ with~}~ 
B_j = (-1)^{j-1}b_{2j}.$$
Applying the inversion relation of Lemma
\ref{lem:inversion} to (\ref{eq:alpha}), one
obtains
\begin{equation}
\alpha_j=(-1)^j\sum_{k=0}^j(-1)^k{2j+2\choose 2k}(2^{2k}-2)B_k.
\label{eq:solalpha}
\end{equation}
On the other hand we have the following relation
for the Bernoulli numbers,
\begin{lemma}(see {\rm \cite{W}})
Denote $
v_j=2(2^{2j}-1)B_j$. 
We have the following recurrence relation :
\begin{displaymath}
v_k-{1\over 2}{2k\choose 2}v_{k-1}+{1\over 2}{2k\choose 4}v_{k-2}-
\cdots+(-1)^{k-1}{1\over 2}{2k\choose 2k-2}
v_1+(-1)^k k=0.
\end{displaymath}
\end{lemma}
Applying the above relation to the
right hand side of (\ref{eq:solalpha}), we have 
\begin{displaymath}
\alpha_j=2(2^{2j+2}-1)B_{j+1}.
\end{displaymath}
Therefore the following result follows
immediately:
\begin{proposition} \label{pro:Yrel}
The following
relations of 
$Y_k$s hold in $\goth{OA}_L$, hence 
in $\goth{OA}$
\begin{equation}
Y_{2n} = \sum_{k \geq n} (-1)^{k-n+1}
{(4^{k-n+1}-1)B_{k-n+1}\over k-n+1}
{2k \choose 2n-1 } Y_{2k+1} \ .\label{eq:even}
\end{equation}
\end{proposition}
Through the morphism (\ref{eq:ta}) for $a=1$, 
$\widehat{\goth{OA}}$ is embedded in the
formal algebra $\goth{sl}_2[[u]]$ where $u$ is
the local coordinate of $t$-plane near $t=1$ 
$$ 
u= t-1 \ . 
$$
Another local coordinate near $t=1$ is given by
$$
v = t^{-1} - 1 
$$
with the relation $\displaystyle{u= \frac{-v}{1+v}, v=
\frac{-u}{1+u}}$. We have $\goth{sl}_2[[u]]
= \goth{sl}_2[[v]]$.  
In the $\goth{sl}_2$-formal algebra, 
the elements $X_k, Y_k$ have the following
symmetric expressions. 
\begin{lemma} \label{lem:XYt}
The generators
$X_k, Y_k (k \geq 0)$ of $\widehat{\goth{OA}}$
can be represented by 
$$
X_k = 2u^k e + 2v^k f \ , \ \ 
\ \  \ Y_k = (-1)^k \left(u^k-
v^k\right) h \
 .
$$
\end{lemma}
{\bf Proof.} The following 
relations hold in $\widehat{\goth{OA}}$ for 
$n \geq 0$
\begin{eqnarray*}
\sum_{k\geq 0} {n \choose k} X_k &= A_n & = 
2 t^n e + 2 t^{-n} f \ , \ 
\cr
\sum_{k\geq 0} {n \choose k}(-1)^k Y_k 
&= G_n & =
(t^n - t^{-n}) h  \ .
\end{eqnarray*}
By an induction procedure, our results follow 
from the identities:
\begin{eqnarray*}
t^n = \left(u +1\right)^n = \sum_{k\geq 0} 
{n \choose k} u^k \ , \\ 
t^{-n} = (v +1)^n = \sum_{k\geq 0}
{ n \choose k } v^k \ .
\end{eqnarray*}
$\Box$ \par \vspace{.2in} \noindent
Using Lemma \ref{lem:XYt} and Proposition
\ref{pro:Yrel}, one obtains the following result.
\begin{proposition}\label{pro:dOAL}
$\goth{OA}_L $ is a
solvable Lie algebra of dimension $L+ [{L\over 2}]$ 
with a basis consisting of $X_k$ and $ Y_j$ with $
0 \leq k , j < L,$ and $j$: odd.
\end{proposition}
Let us examine the structure of the quotient 
$\goth{OA}_{2l}$ more 
closely.
The quotient $\goth{OA}_{2l}/[\goth{OA}_{2l},\goth{OA}_{2l}]$ is
one dimensional and is spanned by $X_0$.
Following the recipe of Malcev \cite{M} for studying the structure of
solvable Lie algebras, we examine the spectrum of
the operator 
${\rm ad}X_0$.
The following result was found by looking in
details at  the
quotient algebras
$\goth{OA}_{2l}$ for small $l$ and then was proved by direct calculation.
\begin{lemma}
The operator ${\rm ad}X_0$ has the eigenvalues 
$0$, $\pm4$ on the space 
$\goth{OA}_{2l}$.
Each of the eigenspace is $l$-dimensional.
A basis of each eigenspace 
is given as follows:
\begin{eqnarray*}
0&:&X_0,\sum_{k=2j-2}^{2l-1}(-1)^k{k-1\choose 2j-4}X_k, \quad 
2\leq j \leq l,\cr
\pm 4&:& 2Y_{2j-1}\pm\left(X_{2j-1}-\sum_{k=2j-1}^{2l-1}
(-1)^k{k-1\choose 2j-2}X_k\right), 
\quad 1\leq j\leq l.
\end{eqnarray*}
Further the 0-eigenvectors 
commute each other.
\end{lemma}
Unfortunately the other commutation relations
among above eigenvectors are not so easy to
determine. Therefore we proceed in the following
manner. Set
\begin{eqnarray}
&&H_0={1\over 2}X_0, \nonumber \\
&&E_0={1\over 4}Y_1+{1\over 8}\left(X_1-\sum_{k=1}^{2l-1}(-1)^kX_k\right),
\label{eq:e0}  \\
&&F_0={1\over 4}Y_1-{1\over 8}\left(X_1-\sum_{k=1}^{2l-1}(-1)^kX_k\right).\nonumber 
\end{eqnarray}
Then they satisfy the relations
$$
[H_0,E_0]=2E_0,\quad [H_0,F_0]=-2F_0.
$$
Set
$$
H_1=\left[E_0,F_0\right] \ , 
$$
and define inductively
\begin{equation}
E_{j+1}={1\over 2}[H_1,E_j],\quad F_{j+1}={-1\over
2}[H_1,F_j]  ,
\quad 0\leq j\leq l-2.\label{eq:induction}
\end{equation}
Then we have
\begin{lemma}\label{lem:indep}
$
\left[E_{j+1},F_{k-1}\right]=\left[E_j,F_k\right].
$
\end{lemma}
{\bf Proof.}
Using the definition of $E_{j+1}$, we have
\begin{eqnarray*}
&[E_{j+1},F_{k-1}]&=[{1\over 2}[H_1,E_j],F_{k-1}] \cr
 &&
=-{1\over 2}[[E_j,F_{k-1}],H_1]-
{1\over 2}[[F_{k-1},H_1],E_j].
\end{eqnarray*}
Since $[E_j,F_{k-1}]$ belongs to $0$-eigenspace
of ${\rm ad}X_0$, the right hand side in the
above becomes 
$$
{-1\over 2}[E_j,[H_1,F_{k-1}]]=[E_j,F_k].
$$
$\Box$ \par \vspace{.2in} \noindent
By Lemma \ref{lem:indep}, the commutator
$[E_j,F_k]$ depends only on
$j+k$, which we denote by
$H_{j+k+1}:= [E_j,F_k]$.
By Proposition \ref{pro:dOAL}, 
a basis of $\goth{OA}_{2l}$ is given by
$X_j$ $0\le j\le 2l-1$ and $Y_{2k+1}$ $0\le k<l$.
Through their relations with $X_j, Y_j$,
one can see that  
$H_j, E_j, F_j$,
$0\leq j\leq l-1$ also constitute a basis of 
$\goth{OA}_{2l}$. 
Further in view of the eigenvalue distribution of ${\rm ad}X_0$, it is easy to
see that $E_j$s and $F_j$s commute among themselves respectively.
Summarizing (\ref{eq:e0}),
 (\ref{eq:induction})
and Lemma \ref{lem:indep} , we have
\begin{theorem} \label{th:EFH}
The elements $E_j$, $F_j$, $H_j$, $0\leq j\leq
l-1$ form a basis of the quotient algebra
$\goth{OA}_{2l}$, in which the following
commutation relations hold,
\begin{eqnarray*}
&&[E_j,F_k]=H_{j+k+1}, \quad [H_j,E_k]=2E_{j+k}, \quad [H_j,F_k]=-2F_{j+k},\cr
&&[E\j,E_k]=0, \quad [F_j,F_k]=0,
\end{eqnarray*}
here if the indices exceed $l-1$ then the corresponding elements are 
regarded to be 0. 
\end{theorem}

The above commutation relations suggests that the
structure of
$\goth{OA}_{2l}$
is very close to $\left({\bf C}[x]/x^l{\bf C}[x]\right)\otimes\goth{sl}_2$.
However the actual structure differs slightly. 
To fix this we
consider the formal algebra
$\widehat{\goth{OA}}$ and employ the loop 
representation of $\goth{OA}$. 
Instead of the variables, $u = t-1 ,  v =
t^{-1}-1$, which we used before,
a convenient coordinate system near $t=1$ for our
purpose now is the following one 
$$
\lambda ={t-t^{-1}\over 2} \ . 
$$
One has $\CZ[[u]] = \CZ[[v]] = \CZ[[ \lambda ]]
$. In fact, the relations of $\lambda$ and
$u, v$  are given by
\begin{eqnarray*}
&&u  = \lambda -1 +
\sqrt{1+\lambda^2} \ , 
\\ &&v  = 
- \lambda -1 + \sqrt{1+\lambda^2} \ , 
\\  &&\lambda
 = {u(2+u)\over 2(1+u)} = {-v(2+v)\over 2(1+v)}
\ .
\end{eqnarray*}
Using Lemma \ref{lem:XYt} and 
substituting the above relations into 
(\ref{eq:e0}) ,
we have
\begin{eqnarray*}
 &&H_0=e+f,\cr
&&E_0={\lambda\over 2}\left(-h+(e-f)\right),\cr
&&F_0={\lambda\over 2}\left(-h-(e-f)\right).
\end{eqnarray*}
We introduce the following
elements in $\goth{sl}_2$ (which correspond
to Pauli matrices in the canonical 
representation of $\goth{sl}_2$)
\bea(l)
\sigma^1= e+f \ , \ \ \sigma^2= -\sqrt{-1}e+
\sqrt{-1}  f \ , \ \ \sigma^3 = h \ .
\elea(sigma)
By the definition of $E_j$, $F_j$ and
$H_j$, we find that
\begin{eqnarray*}
&&H_j=\lambda^{2j}\sigma^1 ,\cr
&&E_j=
\lambda^{2j+1}{\sqrt{-1}\sigma^2 - \sigma^3 \over
2},
\cr 
&&F_j=\lambda^{2j+1}
{-\sqrt{-1}\sigma^2 - \sigma^3 \over 2}.
\end{eqnarray*}
Using the automorphism of $\goth{sl}_2$ by cyclic 
permuting $\sigma^j$s, we obtain an isomorphism
which gives the structure of $\goth{OA}_{2l}$,
\begin{equation}
\goth{OA}_{2l}
\simeq
\bigoplus_{j=0}^{l-1}{\bf C}\lambda^{2j}h+
\bigoplus_{j=0}^{l-1}{\bf C}\lambda^{2j+1}e+
\bigoplus_{j=0}^{l-1}{\bf
C}\lambda^{2j+1}f.\label{eq:homm}
\end{equation}
Note that the above $\goth{OA}_{2l}$ is a
solvable Lie algebra of dimension $3l$. The
derived ideal $[\goth{OA}_{2l},\goth{OA}_{2l}]$
is a nilpotent ideal of dimension $3l-1$.
The classification of nilpotent Lie
algebras is, in general, known to be a  wild
problem.  
Santharoubane has proposed a program of 
classifying nilpotent Lie algebras in \cite{S}.
By analyzing the commutation relations of root vectors of nilpotent
Lie algebras, Santharoubane associates a generalized Cartan matrix (GCM)
to each nilpotent Lie algebra and reduced the classificaton problem
of nilpotent Lie algebras to the one of certain
ideals in the nilpotent part of Kac-Moody
algebras.
However this problem is not easy though.
Santharoubane and others try to classify these
ideals in classical Lie algebras.
For affine Lie algebras the work is done for
$A_1^{(1)}$ and
$A_2^{(2)}$.
According to his classification there are three
series of nilpotent Lie  algebras associated to
$A_1^{(1)}$. One of these series denoted by
${\cal A}_{1,l-1,1}$ is isomorphic to 
\begin{equation}
{\bf C}e+\bigoplus_{j=1}^{l-1}x^j\goth{sl}_2+{\bf C}x^lf.\label{eq:nill}
\end{equation}
We can give an explicit isomorphism
between this Lie algebra and the derived ideal of
$\goth{OA}_{2l}$. However since the formula is
quite complicated, we will not give it here.
Recall here that the affine Lie algebra $A_1^{(1)}$ has several
realizations. The 
most well-known one is the homegeneous realization
\begin{displaymath}
A_1^{(1)}\simeq 
({\bf C}[x]\otimes\goth{sl}_2) \oplus{\bf C}c.
\end{displaymath}
There also exists the principal realization
\begin{displaymath}
A_1^{(1)}\simeq 
({\bf C}[y^2,y^{-2}]\otimes h ) \oplus 
( y{\bf
C}[y^2,y^{-2}]
\otimes e ) \oplus 
( y{\bf C}[y^2,y^{-2}]\otimes f ) .
\end{displaymath}
The presentation of nilpotent Lie algebra (\ref{eq:nill})
refers to homegeneous realization; while the
nilpotent Lie algebra appearing as the derived
ideal of (\ref{eq:homm}) refers to the
principal  realization of $A_1^{(1)}$. In
conclusion we found that the series of nilpotent
Lie algebras 
${\cal A}_{1,l,1}$ appears as the derived ideal of the quotient of
the Onsager algebra $\goth{OA}_{2l}$.

Now we go back to the structure problem of
a general $\goth{OA}_{L}$ and the formal algebra 
$\widehat{\goth{OA}}$. In the following, we
present a slightly different approach from the
previous one by considering more on the loop
structure of 
$\goth{OA}$. First we need the following relations
of powers of the coordinates, $ u, v, \lambda$.
\begin{lemma} \label{lem:uvlam} 
In the Laurent series ring $\CZ((\lambda))$, for
$ n \in \ZZ $, we have
\begin{eqnarray*} 
&&  u^{2n} + v^{2n}, u^{2n-1} + v^{2n-1}  \in
\lambda^{2n}\CZ[[\lambda^2]] , \\
&&  u^{2n} - v^{2n} , u^{2n+1} - v^{2n+1}  \in
\lambda^{2n+1}\CZ[[\lambda^2]] \ .
\end{eqnarray*} 
In fact, the
following ratios tend to 1 as $\lambda
\rightarrow 0 $ 
\begin{eqnarray*}
\frac{u^{2n} + v^{2n}}{2 \lambda^{2m}} \ , \ \ 
\frac{u^{2n-1} + v^{2n-1} }{(2n-1)\lambda^{2n}}
 \ ,
\frac{u^{2n} - v^{2n}}{2n \lambda^{2n+1}} \  \ 
,
\ \ 
\frac{u^{2n+1} - v^{2n+1}}{2 \lambda^{2n+1}} \ .
\end{eqnarray*}
\end{lemma}
\noindent Now we derive the structure of
$\widehat{\goth{OA}}, 
\goth{OA}_L$  as follows:
\begin{theorem}
Denote  the Lie
subalgebras of
$ \goth{sl}_2[[\lambda]]$ for $L \geq 0$ 
\begin{eqnarray*}
&&\goth{sl}_2 \langle\langle \lambda
\rangle\rangle  : = 
\CZ[[\lambda^2 ]] h +  
\lambda \CZ[[\lambda^2 ]] e + \lambda \CZ[[\lambda^2 ]]
f \ \subset \goth{sl}_2[[\lambda]] , \\ &&
\goth{sl}_2\langle\langle \lambda \rangle\rangle^L : = 
\goth{sl}_2\langle\langle \lambda \rangle\rangle \bigcap \lambda^L
\goth{sl}_2[[\lambda]] \ .
\end{eqnarray*} 
Then we have the isomorphism 
$$
\widehat{\goth{ OA}} \ \ \simeq  \ \ 
\goth{sl}_2 \langle\langle \lambda
\rangle\rangle \ ,
$$
which induces the isomorphisms 
$$
\begin{array}{ll}
\widehat{\goth{ OA}}^L \ \ \simeq \ \  
\goth{sl}_2 \langle\langle \lambda
\rangle\rangle^L \ , \ \ 
& \goth{ OA}_L   \ 
 \simeq \ 
\goth{sl}_2 \langle\langle \lambda
\rangle\rangle/\goth{sl}_2 \langle\langle \lambda
\rangle\rangle^L \ .
\end{array}
$$
\end{theorem}
{\bf Proof.} With the elements $\sigma^j$s
(\req(sigma)) of $\goth{sl}_2$ as before, there
is a natural isomorphim 
$$
\goth{sl}_2 \langle\langle \lambda
\rangle\rangle 
\ \ \simeq \ \  
\goth{sl}_2 \langle\langle \lambda
\rangle\rangle':= 
\CZ[[\lambda^2]]
\sigma^1 +
\lambda \CZ[[\lambda^2]]
\sigma^2+ \lambda \CZ[[\lambda^2]]\sigma^3 \ .
$$
We need only to show the results by replacing 
$\goth{sl}_2 \langle\langle \lambda
\rangle\rangle $ as 
$\goth{sl}_2 \langle\langle \lambda
\rangle\rangle' $. The expressions of 
$X_k, Y_k$ in Lemma \ref{lem:XYt} becomes
$$
X_k =  (u^k +v^k)\sigma^1 + \sqrt{-1}(u^k -v^k)
\sigma^2   \ , \ \ 
\ \  \ Y_k = \ (-1)^k (u^k-
v^k) \sigma^3 \ , 
$$
by which and relations in 
 Lemma \ref{lem:uvlam}, $\widehat{\goth{OA}}$ is a
subalgebra of
$\goth{sl}_2 \langle\langle \lambda
\rangle\rangle'$ and 
$\widehat{\goth{OA}}^L  \subseteq \goth{sl}_2
\langle\langle \lambda \rangle\rangle'^L$. There
is an induced canonical morphism 
$$
\rho_L:
\goth{OA}_L
= \widehat{\goth{OA}}/\widehat{\goth{OA}}^L
\longrightarrow  \goth{sl}_2 \langle\langle
\lambda \rangle\rangle'/\goth{sl}_2
\langle\langle \lambda
\rangle\rangle'^L \ 
$$
for $ L \geq 0 $. 
It remains to show that $\rho_L$ is an isomorphism. By
Proposition \ref{pro:dOAL},  
$ X_k,  Y_{2j+1}, (0 \leq k, 2j+1 <L)$ form a
basis of  
$\goth{OA}_L$. 
While in 
$\goth{sl}_2 \langle\langle \lambda
\rangle\rangle'/\goth{sl}_2 \langle\langle \lambda
\rangle\rangle'^L$,  
$\lambda^{2l} \sigma^1, \lambda^{2j+1} \sigma^2,
\lambda^{2j+1}\sigma^3$, $ 0 \leq 2l, 2j+1 <L$
form a basis. By using the behavior of ratios in 
 Lemma \ref{lem:uvlam}  near
$\lambda=1$, the matrix of $\rho_L$ for these
bases is an invertible lower
triangular one. Hence $\rho_L$ is an isomorphism.
$\Box$ \par \vspace{.2in} \noindent
{\bf Remark.} The isomorphism in the above
theorem for 
$\goth{OA}_L$ for $ L=2l $ is the one in 
(\ref{eq:homm}). For $L=2l+1$, it is given by
\begin{eqnarray*}
\goth{OA}_{2l+1}
\simeq
\bigoplus_{j=0}^{l}{\bf C}\lambda^{2j}h+
\bigoplus_{j=0}^{l-1}{\bf C}\lambda^{2j+1}e+
\bigoplus_{j=0}^{l-1}{\bf
C}\lambda^{2j+1}f \ ,
\end{eqnarray*}
which can be derived by the same argument as
before for  (\ref{eq:homm}) by examining the
eigenvectors of ad$X_0$.

\section{Irreducible Representations
of the Onsager Algebra and Superintegrable Chiral
Potts Model} 
In this section, we are going
to derive the  classification of irreducible
representations of the Onsager algebra. For ${\bf
a} = (a_1,
\ldots, a_n) \in \CZ^{*n}$, one has the
evaluation   morphism of $\goth{OA}$ into the sum
of $n$ copies of $\goth{sl}_2$ defined by 
$$
ev_{\bf a} : \goth{OA} \longrightarrow \bigoplus^n \goth{sl}_2
\ , \ \ 
\ \ X \mapsto ( ev_{a_1}(X), \ldots, ev_{a_n}(X)) \ 
,
$$
where $ev_{a_j}$ is the evaluation of
$\goth{OA}$ at $a_j$.  Denote 
$$
U_{\bf a}(t):= \prod_{a \in \{a_1, \ldots, a_n \} }
U_a(t) 
\in \CZ[t] \ .
$$ 
\begin{lemma}\label{lem:param}
We have ${\rm Ker}(ev_{\bf a}) = \goth{I}_{U_{\bf
a}(t)}$. The surjectivity of 
$ev_{\bf a}$ is
equivalent to  $a_j 
\neq \pm1 , a_j \neq a_k^{\pm1}$ for $j \neq k $.  
\end{lemma}
{\bf Proof.} For the determination of ${\rm
Ker}(ev_{\bf a})$,
through the diagonal map $ \triangle :\goth{sl}_2
\longrightarrow  \goth{sl}_2 \oplus \goth{sl}_2$ and 
the involution $\theta : \goth{sl}_2
\longrightarrow  \goth{sl}_2$ , one may assume that
 ${\bf a}=(a_1, \ldots, a_n)$ satisfies the 
condition, $
 a_j \neq
a_k^{\pm1}$ for $ j \neq k$. 
In fact, for simplicity let us take  $n=2$ as an
example; a similar argument can apply to the  
case of a general $n$. When
$a_1=a_2= a$, the evaluation map 
$ev_{(a,a)}$ of $\goth{OA}$  can
be  reduced to $ev_a$ by the relation $ ev_{(a,a)} =
\triangle ev_a$. When $a_1 = a_2^{-1} =a$, one 
reduces  the map $ev_{(a, a^{-1})}$ to $ev_{(a,a)}$
by 
$ev_{(a, a^{-1})} = ({\rm id}, \theta) ev_{(a,a)}$. Hence
both situations  are reduced to the case
$n=1$.  By  Lemma \ref{lem:ker},
the argument also shows the non-surjectivity of $ev_{\bf
a}$ if   
${\bf a}$ has two
components with equal or reciprocal  values; the
same conclusion for ${\bf a}$ with one component
equal to $\pm1$. Conversely, by Lemma
\ref{lem:ker} $ev_{\bf a}$ is  surjective for
${\bf a}$ with $a_j \neq \pm1$ and $
 a_j \neq
a_k^{\pm1}$ for $j \neq k$. 
Now we may assume
${\bf a}$ with $
 a_j \neq
a_k^{\pm1}$. Then 
$$
{\rm Ker}(ev_{\bf a}) = \bigcap_{j=1}^n {\rm Ker}(ev_{a_j}) 
= \bigcap_{j=1}^n  \goth{I}_{U_{a_j}(t)} = 
\goth{I}_{U_{\bf a }(t)} \ .
$$
$\Box$ \par \vspace{.2in} \noindent 
The quotient space of
$\CZ^*$ by identifying $a$ with $a^{-1}$ is again 
parametrized by $\CZ^*$ with the variable denoted
by 
$\hat{a}
\in \CZ^*$, which is related to $a$ by the
rational map  
$$
\CZ^* \longrightarrow \CZ^* \ , \ \ a \mapsto \hat{a}: = 
\frac{a+a^{-1}}{2} \ .
$$ 
By composing with the involution $\iota$ of 
$\goth{OA}$, the representations
$ev_{a^{-1}}$ and $ ev_a$ are equivalent   
$ev_{a^{-1}}= ev_{a}
\iota$. 
We shall use the following  convention if no confusion
could arise:
$$
\epsilon_{\hat{a}}: = ev_{a} \ ; \ \ \ 
\epsilon_{\hat{\bf a}}: = ( \epsilon_{\hat{a_1}}, 
\ldots , \epsilon_{\hat{a_n}} ) \ , \ \ {\rm for} \ 
\hat{\bf a}= (\hat{a_1}, 
\ldots , \hat{a_n} )
$$
Denote ${\cal S}$ the collection of all
the non-trivial integral representations of
$\goth{sl}_2$. It is known that elements in
${\cal S}$  are  labelled by positive
half-integers $s$, which corresponds to the
irreducible representation of $\goth{sl}_2$  on
the 
$(2s+1)$-dimensional vector space $V(s)$.
Therefore the  effective  irreducible integral
representations of
$\stackrel{n}{\bigoplus} \goth{sl}_2$ are indexed by  ${\bf
s}=(s_1,
\ldots , s_n)
\in {\cal S}^n $, where  $\stackrel{n}{\bigoplus}
\goth{sl}_2$ acts on the vector space $V({\bf s})
(:= 
\otimes V(s_j) )$  by the relation,
$$
(x_1, \ldots, x_n) v = \sum_j ( 1 \otimes \ldots
\otimes  x_j \otimes \ldots \otimes 1) v \ , \ \ 
x_j \in \goth{sl}_2 \ , \ v \in V({\bf s}) \ .
$$ 
Combining the above representation with $
\epsilon_{\hat{\bf a}}$, one obtains a representation of 
$\goth{OA}$ on $V({\bf s})$ 
$$
\rho_{(\hat{\bf a }, {\bf s})} : 
\goth{OA} \longrightarrow \goth{gl}(V({\bf s})) \
,
\
\ (\hat{\bf a }, {\bf s}) \in \CZ^{*n} \times 
{\cal S}^n \ . 
$$ 
The Hermitian condition of $\rho_{(\hat{\bf a }, {\bf
s})}$ is given by  $|a_j|=1$, equivalently, 
$\hat{a_j}$ in the real interval $[-1, 1]$ for
all $j$, {\it i.e.}  $
{\bf a}=  ( e^{\sqrt{-1}\theta_1}, \ldots,e^{\sqrt{-1}}\theta_n),
 \ \hat{\bf a}=   ( \cos
(\theta_1), \ldots, \cos(\theta_n))$.  Denote 
$$
\begin{array}{ll}
{\cal C}_n := & \{ \hat{\bf a } = ( \hat{a}_1, \ldots, 
\hat{a}_n ) \in (\CZ^*
\setminus
\{\pm 1\})^n \ | \ \hat{a}_j \neq \hat{a}_k \ \ {\rm for} \
j \neq k \} \ , \\
{\cal D}_n :=& {\cal C}_n \bigcap (-1, 1)^n \ .
\end{array}
$$
By Theorem 2 and the structure of ${\rm
Ker}(ev_{\bf a})$  in Lemma \ref{lem:param}, one 
obtains the following  results:
\begin{proposition} \ \label{prop:Onsrep}

(I).  ${\rm Ker}(\rho_{(\hat{\bf a }, {\bf s})}) = \goth{I}_{U_{\bf a}(t)}$ 
for $(\hat{\bf a }, {\bf s}) \in \CZ^{*n} 
\times  {\cal S}^n$.

(II). $\rho_{(\hat{\bf a }, {\bf s})}$
is irreducible  if and only if $\hat{\bf a }
\in {\cal C}_n$. 

(III). $\rho_{(\hat{\bf a }, {\bf s})}$
is irreducible Hermitian if and only if $\hat{\bf a } 
\in {\cal D}_n$. 

(IV). For $(\hat{\bf a }, {\bf s})
\in {\cal C}_n \times {\cal
S}^n, (\hat{\bf a'},
{\bf s'}) \in {\cal C}_{n'} \times  {\cal
S}^{n'} $, the representations 
$\rho_{(\hat{\bf a}, {\bf s})}, \rho_{(\hat{\bf a' }, {\bf
s'})}$ are equivalent if and only if $ n=n' $ and $
\hat{a}_j'= \hat{a}_{\sigma(j)}, s_j'= s_{\sigma(j)} $
for some permutation $\sigma$ of indices.
\end{proposition} 
Now  we are going to 
classify  all the irreducible representations
of $\goth{OA}$.
\begin{lemma}\label{lem:irred}
Let $\rho$ be a non-trivial
irreducible representation  of $\goth{OA}$ in 
$\goth{sl}(V)$. Then 
${\rm Ker}(\rho) = \goth{I}_{U_{\bf a}(t)}$ for some  
${\bf a} = (a_1, \ldots, a_n) 
\in (\CZ^* \setminus \{\pm1\})^n$ with $a_j \neq
a_k^{\pm1}$  for $j \neq k$.
\end{lemma}
{\bf Proof.} By 
Schur's Lemma, ${\rm Ker}(\rho)$ is a closed
ideal.  By Theorem \ref{th:central}, ${\rm
Ker}(\rho)=
\goth{I}_{P(t)}$ for some 
$ P(t) = \prod_{j=1}^n
U_{a_j}(t)^{m_j}$, where $ a_j \in
\CZ^*, m_j
\in \ZZ_{>0}$ with $a_j \neq a_k^{\pm 1} $ for $ j \neq
k$, and $m_j$ even whenever $a_j = \pm1$.
It  suffices to show $m_j= 1$ for all $j$. 
Otherwise, we
may  assume $m_1 \geq 2$. Define the 
polynomial $R(t)=U_{a_1}(t)^{m_1-1}\prod_{j>1}U_{a_j}(t)^{m_j}$. 
Then  $\goth{I}_{R(t)}/\goth{I}_{P(t)}$ is a
non-trivial abelian ideal of 
$\goth{OA}/\goth{I}_{P(t)}$. By the
irreducibility of $\rho$, one has $V=
(\goth{I}_{R(t)}/\goth{I}_{P(t)})V$.   Let $V=
V_1 \oplus
\ldots 
\oplus V_r$ be the eigenspace decomposition of $V$ 
with respect to the
$\goth{I}_{R(t)}/\goth{I}_{P(t)}$-action, then
eigenvalues
$\lambda_j$ for $V_j 
(1 \leq j \leq r)$ are distinct. Here
$\lambda_j$ is  a linear functional on
$\goth{I}_{R(t)}/\goth{I}_{P(t)}$. As the
representation takes value in 
$\goth{sl}(V)$,  the number $r$ is at least $2$.
We are going to show that the vector space $V_1$
gives rise to a subrepresentation of
$\goth{OA}$, hence a contradiction to the
irreduciblity of $\rho$. Let $v$ be en element of 
$V_1$ and $X
\in
\goth{OA}$.  Denote
 $
\rho(X)(v) = \sum_{l=1}^r v_l$ with $ v_j \in V_j
$. 
For each $j \geq 2$, we choose an element $Z_j 
\in \goth{I}_{R(t)}$ such that the class of $Z_j$
in $\goth{I}_{R(t)}/\goth{I}_{P(t)}$ takes
different values for $\lambda_1$ and $\lambda_j$,
$\lambda_1( Z_j + \goth{I}_{P(t)}) \neq 
 \lambda_j( Z_j + \goth{I}_{P(t)})$. 
As $\rho(X) \rho(Z_j) - \rho(Z_j) \rho(X) ( = 
\rho([X, Z_j]) $ is an element of 
$\rho(\goth{I}_{R(t)})$,  we have 
$$
V_1 \ni \rho(X) \rho(Z_j)(v) - \rho(Z_j) \rho(X)(v) = 
\lambda_1( Z_j + \goth{I}_{P(t)}) \sum_{l=1}^r v_l - 
\sum_{l=1}^r \lambda_l( Z_j + \goth{I}_{P(t)})v_l \ , 
$$
which implies $v_j=0$ for $j \geq 2$. Therefore
$V_1$ is a representation of $\goth{OA}$.   
$\Box$ \par \vspace{.2in} \noindent 
Now we can derive the following result in
\cite{D90, R91}.
\begin{theorem} \label{th:irr}
Any non-trivial
irreducible representation  of $\goth{OA}$ 
is represented by $
\rho_{(\hat{\bf a }, {\bf s})}$ for some $( \hat{\bf a
}, {\bf s}) \in  {\cal C}_n \times {\cal S}^n $ , $n \in
\ZZ_{>0}$. Subsequently,  all the
irreducible Hermitian
 representations  of $\goth{OA}$ are given by $
\rho_{(\hat{\bf a }, {\bf s})}$ for $( \hat{\bf a },
{\bf s}) \in
\bigsqcup_{n \in
\ZZ_{>0}} ( {\cal D}_n \times {\cal S}^n) $, modulo the
following relation:
$$
( \hat{\bf a }, {\bf s})= ((
\hat{ a}_1, \ldots, \hat{a}_n)), (s_1, \ldots, s_n) )
\sim ( \hat{\bf a }', {\bf s}') =
((
\hat{ a}'_1, \ldots, \hat{a}'_n)), (s'_1, \ldots, s'_n) ) \
, $$
where
$
\hat{a}_j'= \hat{a}_{\sigma(j)}, s_j'= s_{\sigma(j)} \ ,
$ 
and $\sigma$ is a permutation of indices.
\end{theorem}
{\bf Proof.} For a non-trivial irreducible
representation $\rho$ of $\goth{OA}$, by Lemma
\ref{lem:prod},\ref{lem:ker} and
\ref{lem:irred},  the Lie-algebra
$\goth{OA}/{\rm Ker}(\rho)$  is isomorphic to
$\stackrel{n}{\oplus}\goth{sl}_2$ for some
positive integer $n$. As an irreducible
representation  of
$\stackrel{n}{\oplus}\goth{sl}_2$ is obtained by
tensoring irreducible one of its facors,
the result  follows immediately.
$\Box$ \par \vspace{.2in} \noindent 
{\bf Remark.} As it is known in the affine Lie
algebra theory, all irreducible finite
dimensional represenations of $\goth{sl}_2$-loop
algebra ( or equivalently $\goth{sl}_2$-affine
algebra ) 
are isomorphic to tensor products of irreducible ones of
$\goth{sl}_2$ through the 
evaluation at distinct non-zero complex numbers
$a_j$s (see e.g., \cite{Ka}  Exercise 12. 11 ). 
While for the irreucible representations of
Onsager algebra $\goth{OA}$, one requires all the 
values $a_j \neq \pm 1$, plus the identification
of $a_j$ and $a_j^{-1}$ which produce the same
representation of
$\goth{OA}$.  Hence the discussion of Sect. 5 on
the effect when the Onsager algebra  is valued
at $\pm 1$  reveals the essence of the Onsager
algebra different from $\goth{sl}_2$-loop
algebra from the representation point of view.

We now discuss a physical application of
the previous results to superintegrable chiral
 Potts $N$-state model. The Hamiltonian is the
spin chain  of a  
parameter $k'$ \cite{AMPT}-\cite{BBP},
\cite{GR, P},
$$
H(k') = H_0 + k' H_1 \ , 
$$
with $H_0, H_1$ the 
Hermitian operators
acting on the vector space of 
$L$-tensor of
$\CZ^N$, defined by  
$$
H_0 = -2 \sum_{l=1}^L \sum_{n=1}^{N-1} (1-\omega^{-n})^{-1}
X_l^n , \ \  \  
 H_1 = -2 \sum_{l=1}^L \sum_{n=1}^{N-1} (1-\omega^{-n})^{-1}
Z_l^nZ_{l+1}^{N-n} \ ,
$$
where $\omega = e^{\frac{2 \pi
\sqrt{-1}}{N}} $,  
$
X_l = I \otimes \ldots \otimes \stackrel{l{\rm
th}}{ X}
\otimes \ldots \otimes I $, $ 
Z_l = I \otimes \ldots \otimes \stackrel{l {\rm
th}}{ Z}
\otimes \ldots \otimes I ,  (Z_{L+1}= Z_1 ) $.  
Here $I$ is the identity operator,   
$X, Z$ are the operators of $\CZ^N$ with 
the relation, $ZX = \omega XZ$, 
defined by  $X|m\rangle= |m+1\rangle ,
Z|m\rangle = \omega^m |m\rangle$,  
$m \in \ZZ_N$. The operator
$H(k')$ is  Hermitian for real $k'$, hence
with the real eigenvalues. 
It is 
the Ising quantum chain \cite{K} for $N=2$.
For $N=3$, one obtains the 
$\ZZ_3$-symmetrical self-dual chiral clock model 
with the chiral angles $
\frac{\pi }{2}$, 
which was studied by  
Howes, Kadanoff and M. den
Nijs \cite{HKN}. 
For a general $N$, we set 
$$ 
A_0 = -2 N^{-1} H_0 , \ \  A_1
= -2 N^{-1} H_0 \ ,
$$
then $A_0, A_1$ satisfy the $DG$-condition, which
by Theorem \ref{th:irr} and Proposition
\ref{prop:Onsrep},  ensures that 
the eigenvalues of the unitary
operator $H(k')$ have  the following special form
as in the Ising model,
$$
 a  + b k'   + 2N
\sum_{j=1}^{n} m_j 
\sqrt{1+k'^2 - 2k' \cos (\theta_j) }
\ ,
$$ 
where $a, b, \theta_j \in \RZ$, and 
$m_j =  -s_j, (-s_j +1) , \ldots ,  s_j$ 
with $s_j$ a positive half-integer 
\cite{D90}. 
For the ground-state
sectors, by the computation  of
superintegrable chiral Potts model 
 Baxter \cite{B} obtained the corresponding
eigenvalue of
$H(k')$ given by the spin $\frac{1}{2}$
representation of
$\goth{sl}_2$, i.e., $s_j=\frac{1}{2}$ for all
$j$, with an explicit formula of $\theta_j$s
(see, e.g., \cite{R91}). However these results are
not obvious from the representation theoretic
point of view. The understanding of the exact
form of of eigenvalues of $H(k')$ has still been left as a
theoretical challenge of the study of Onsager
algebra.

\section{ Further Remarks}
In this paper, we have obtained the structure of 
closed ideals $I$ of the Onsager algebra
$\goth{OA}$ and established their relation with
reciprocal polynomials
$P(t)$, $I=
\goth{I}_{P(t)}$. For $P(t)= U_a(t)^L$, we have 
determined the Lie algebra structure of the 
quotient algebra $\goth{OA}/\goth{I}_{P(t)}$.
These 
results, together with a polynomial $P(t)$ of 
mixed types, 
 should make some significant extensions of
our knowledge of    solvable or nilpotent
alegbras. 
Generalizations of the Onsager algebra to other
loop  groups or Kac-Moody algebras, like the one
in
\cite{AS, UI}, should provide ample examples
of solvable algebras. The intimate relation of
Onsager algebra and the superintegrable chiral
Potts model described in Sect.6 also suggests the
potential links of Onsager algebra or its
generalized ones to other quantum integrable
systems in statistical mechanics.
We hope that the further development of 
the subject  will eventually lead to some 
interesting results in Lie-theory with possible 
applications in quantum integrable models.

\section{Acknowledgements}
We would like to thank B. M. McCoy for
numerous illuminating discussions on the relation
of the Onsager algebra and current algebras. We
would also like to thank  
A. N. Kirillov for his interest and encouragement 
while we were writing this paper.

\end{document}